\documentstyle[12pt,amssymb] {amsart}

\begin {document}

\title[Non-Commutative Analytic Toeplitz Algebras]{Factoring in Non-Commutative Analytic Toeplitz Algebras}
\author[D. W. Kribs]{David W. Kribs}
\thanks{Partially supported by an NSERC graduate scholarship}
\address{Pure Math.\ Dept.\\U. Waterloo\\Waterloo, Ont. N2L--3G1\\CANADA}
\email{dwkribs@@barrow.uwaterloo.ca}

\newcommand{\WOT}{\textsc{wot}}
\newcommand{\hinf}{H^{\infty}}
\newcommand{\htwo}{H^{2}}
\newcommand{\hht}{\hat{h}}

\begin{abstract}
The non-commutative analytic Toeplitz algebra is the $\WOT$-closed algebra
generated by the left regular representation of the free semigroup on $n$
generators.
The structure theory of contractions in these algebras is examined. Each is
shown to have an $\hinf$ functional calculus.
The isometries defined by words are shown to factor only as the words do over the
unit ball of the algebra. This turns out to be false over the full algebra.
The natural identification of $\WOT$-closed left ideals with invariant subspaces of the algebra is shown to hold only for a proper subcollection of the subspaces.
\end{abstract}

\maketitle

\theoremstyle{plain}
\newtheorem{thm}{Theorem}[section]
\newtheorem{prop}[thm]{Proposition}
\newtheorem{cor}[thm]{Corollary}
\newtheorem{lem}[thm]{Lemma}

\theoremstyle{definition}

\newtheorem{rem}[thm]{Remark}
\newtheorem{defn}[thm]{Definition}

\newcommand{\I}{{\cal{I}}}
\newcommand{\M}{{\cal{M}}}
\newcommand{\F}{{\cal{F}}}
\renewcommand{\H}{{\cal{H}}}
\newcommand{\T}{{\cal{T}}}
\newcommand{\B}{{\cal{B}}}
\newcommand{\K}{{\cal{K}}}
\newcommand{\gL}{{\frak{L}}}
\newcommand{\gR}{{\frak{R}}}
\newcommand{\bbC}{{\Bbb{C}}}
\newcommand{\bbT}{{\Bbb{T}}}

\newcommand{\czi}{c_{i_{j-1}\cdots i_{k}}}
\newcommand{\Bzi}{B_{i}}
\newcommand{\Bzis}{B_{i_{1}}}
\newcommand{\Pf}{\noindent {\sc Proof.\ }}
\newcommand{\bx}{\hfill $\blacksquare$ \medbreak}
\newcommand{\onef}{\frac{1}{f}}
\newcommand{\Lat}{\operatorname{Lat}}
\newcommand{\Id}{\operatorname{Id}}

In \cite{topinv} and \cite{topalg}, the algebraic and invariant subspace structures of the
non-commutative analytic Toeplitz algebras were developed extensively.
Several analogues of the analytic Toeplitz algebra were obtained. Many
of these results came from a lucid characterization of the $\WOT$-closed
right ideals of these algebras.
Although technical difficulties were encountered, a similar characterization
of the left ideals was expected.
In this paper, it is shown that, although it holds for a subcollection, the analogous characterization of the $\WOT$-closed left
ideals fails.
The reason for this failure is a deep factorization problem in these algebras.
Typically, when norm conditions are placed on possible factors of operators
in these algebras, reasonable factorization results can be obtained. Indeed,
 positive results regarding isometries in the unit ball are
included. However, in the general setting it turns out that even seemingly
obvious unique factorizations do not hold. The examples provided go toward
understanding the fabric of these algebras as well as the pathologies of
factorization involved.
Many of these examples rely on an understanding of the structure theory of
contractions in these algebras. The minimal isometric dilation of these
contractions is determined. Further, each contraction is shown to have an $\hinf$ functional
calculus.
The author would like to thank his supervisor, Ken
Davidson, for all of his assistance.

The terminology and notation used in this paper is that of \cite{topinv} and
\cite{topalg}. The left regular representation of the unital free semigroup $\F_{n}$
on $n$ generators $z_{1},\ldots ,z_{n}$ acts naturally on $\ell_{2}(\F_{n})
= \H_{n}$
by $\lambda\,(w)\,\xi_{v}=\xi_{wv}$, for $v$, $w$ in $\F_{n}$.
The non-commutative analytic Toeplitz algebra $\gL_{n}$ is the unital,
$\WOT$-closed algebra generated by the isometries $L_{i}=\lambda\,(z_{i})$
for $1\leq i \leq n$. For convenience, given $w$ in $\F_{n}$ put $L_{w}=\lambda\,(w)$. The case $n=1$ is exactly the analytic Toeplitz algebra.
Throughout the paper, $n$ is taken to be a finite positive integer. However,
the results of (\cite{topinv},\cite{topalg}) on the structure of $\gL_n$ which
are used here are valid for $n=\infty$. So the results of this paper go through
for $n=\infty$ with only minor notational changes. For ease of presentation,
the paper is written as though $n$ is finite.
The algebra corresponding to the right regular representation is denoted by
$\gR_{n}$. The generating isometries are defined by $\rho(w)=R_{w'}$ where $R_u\xi_v=\xi_{vu}$, and $w'$ denotes the word $w$ in reverse order. It is unitarily equivalent to $\gL_{n}$ and is precisely the
commutant of $\gL_{n}$ (see \cite{topinv} and \cite{po1}).

In the first section it is proved that the minimal isometric dilation of each
non-unitary contraction in $\gL_{n}$ is a shift. This is accomplished by showing that the
powers of the adjoint of such a contraction converge strongly to 0. All of these shifts have infinite multiplicity. Further, each of these contractions has an $\hinf$
functional calculus.

The second section contains positive factorization results for isometries
over the unit ball of $\gL_{n}$. The isometries $L_{w}$ factor exactly as the
words $w$. Hence, such isometries are irreducible over the unit ball exactly
when the length of the word $|w|= 1$.
A characterization of a broader class
of irreducible isometries over the unit ball is also included.
Surprisingly, these isometries are reducible over the full algebra. Particular
factors are constructed using the structure of the orthogonal complement of their ranges.

The last section contains a discussion of the $\WOT$-closed ideals of these algebras.
Unlike right ideals, the left ideal generated by the isometry $L_{2}$ is not
$\WOT$-closed. The same is true for the two-sided ideal generated by $L_{2}$.
The natural identification of subspaces
in $\Lat \gR_{n}$ with the $\WOT$-closed right and two-sided ideals is shown
to only transfer over to a proper subcollection of $\Lat \gL_{n}$, when
considering left ideals.

\section{Contractions}

The inner-outer factorization in $\gL_{n}$ (see \cite{po}, \cite{topinv} and
\cite{po1}) shows that much of the game in these algebras is played in the
world of isometries (inner operators). Thus it is worthwhile to examine the
structure theory of these operators.  It turns out there is much that one can
say about such isometries when they are considered in the more general setting of
contractions.

For notational purposes
 recall from \cite{topalg} that given $k\geq 1$, every $X$ in $\gL_{n}$ can be
 uniquely written as a sum
\begin{eqnarray} X = \sum_{|w|< k}x_{w}L_{w}\,+\,\sum_{|w|= k}L_{w}X_{w},
\end{eqnarray}
 where $x_{w}\in\bbC$ and $X_{w}\in\gL_{n}$ (for $n=\infty$, the first sum
 belongs to $\ell_{2}$ and the second sum is actually a $\WOT$-limit). The
 scalars $\{x_{w}\}_{w\in\F_{n}}$ are called the Fourier coefficients of $X$.
 This notation is justified since they determine the operator. Indeed,
 $X\,\xi_{1}=\sum_{w}x_{w}\xi_{w}$ and hence
 \[ X\,\xi_{v} = X(R_{v}\xi_{1}) = R_{v} (X\xi_{1}) = \sum_{w}x_{w}\xi_{wv}. \]
 One writes $X\sim\sum_{w}x_{w}L_{w}$. For the sake of convenience, given a word
 $w=z_{i_{1}}\cdots z_{i_{k}}$ also let $X_{i_{1}\cdots i_{k}}$ (respectively
$x_{i_{1}\cdots i_{k}}$) denote the operator $X_{w}$ (respectively the scalar $x_{w}$)
 in the expansion (1)
 of $X$. Lastly, the scalar $x_{0}$ is taken to be the inner product $(X\xi_{1},
\xi_{1})$.

The starting point is a simple result which is used several times throughout
the paper. In \cite{topalg} it is shown that the only normal operators
belonging to $\gL_{n}$ are the scalars. Thus the only unitaries in $\gL_{n}$
are scalar. The proof of the latter fact is actually quite elementary.

\begin{prop} \label{unit_prop} The collection of unitary operators in $\gL_{n}$ is exactly the set
$\{\lambda I\, |\, \lambda\in\bbT\}$.
\end{prop}

{\Pf} Let $U$ be unitary in $\gL_{n}$. Then there is a unit vector $\eta$ in
$\H_{n}$ with $U\,\eta = \xi_{1}$. The scalar $\lambda=(\eta ,\xi_{1})$ satisfies
$|\lambda|\leq 1$. Consider the Fourier expansion,
$ U \sim \sum_{w}a_{w}L_{w}.$
Evidently then,
\[ 1 = \| U\,\xi_{1} \|^{2} = \sum_{w}|a_{w}|^{2} \geq |a_{0}|^{2}. \]
However, it is also true that
\[ 1 = \left| (U\eta ,\xi_{1}) \right| = |a_{0}\,\lambda|. \]
Whence $|a_{0}| = |\lambda| = 1$. Thus $\eta=\lambda\,\xi_{1}$, so that           \[ \overline{\lambda}\,\xi_{1} = \overline{\lambda}\,U\,\eta = U\,\xi_{1}. \] Therefore $U = \overline{\lambda}\,I$, and the proof is finished.
{\bx}

In the seminal text \cite{nagfo}, it is shown that understanding the behaviour
of the powers of the adjoint of a contraction is a key issue. In particular,
strong convergence to 0 yields information on the minimal isometric dilation
of the operator, as well as an $\hinf$ functional calculus. This condition
holds for every non-unitary contraction in $\gL_{n}$.

\begin{lem} If $L$ is a non-unitary contraction in $\gL_{n}$,  then
\[ \lim_{k\rightarrow\infty}\| (L^{*})^{k}\xi \| = 0, \]
for all $\xi$ in $\H_{n}$.
\end{lem}

{\Pf} The key is the unique decomposition of $L$. Using (1), write
\[ L = \lambda\, I + \sum_{i=1}^{n}L_{i}A_{i}, \]
with each $A_{i}$ in $\gL_{n}$.
Then by the previous proposition, $| \lambda |<1$ since $L$ is not unitary. Hence if $A=\sum_{i=1}^{n}L_{i}A_{i}$, then $\| A \| < 2$.

The lemma will be proved for basis vectors corresponding to words.
Suppose that $w$ is a word of length $l$.
For $k>l$,
\begin{eqnarray*}
 (L^{*})^{k}\xi_{w}  &=& (\overline{\lambda}\,I + A^{*})^{k}\,\xi_{w} \\                                  &=&  \left( \overline{\lambda}^{k}I + {k \choose 1} \overline{\lambda}^{k-1}A^{*}
  + \ldots + {k \choose  l}
\overline{\lambda}^{k-l}(A^{*})^{l} \right) \,\xi_{w} + 0.
\end{eqnarray*}
For $0\leq j \leq l$ define polynomials $p_{j}$ by
$p_{j}(x) = \frac{1}{j{\rm !}}\, x\,(x-1)\cdots (x-j+1).$
Then one has
\[ (L^{*})^{k}\xi_{w} = \sum_{j=0}^{l}p_{j}(k)\,\overline{\lambda}^{k-j}(
A^{*})^{j}\xi_{w}. \]
However, the limit $\lim_{k\rightarrow\infty}\frac{k^{m}}{\alpha^{k}}=0$  holds
for any real number $m$ and $\alpha > 1$ (see \cite{rudin} p.57). It follows that  for
$0\leq j \leq l$,
\[ \lim_{k\rightarrow\infty}p_{j}(k)\,|\lambda|^{k-j} = 0. \]
Now given $\varepsilon > 0$, choose $K>l$ such that $k\geq K$ implies that
\[ p_{j}(k)\,|\lambda|^{k-j} < \varepsilon \]
for $0\leq j \leq l$. Then for all sufficiently large $k$ one has
\begin{eqnarray*}
\| (L^{*})^{k}\xi_{w} \| &\leq& \sum_{j=0}^{l}p_{j}(k)\,|\lambda|^{k-j}\|
(A^{*})^{j}\xi_{w}\| \\
                         &<& \sum_{j=0}^{l}\varepsilon\,2^{j} \\
                         &=& (2^{l+1}-1)\,\varepsilon .
\end{eqnarray*}
Hence,
 $\lim_{k\rightarrow\infty}\| (L^{*})^{k}\xi_{w} \| =0$ for all words $w$.

That it is true in full generality follows from the boundedness of the
sequence $\{(L^{*})^{k}\}$. Indeed, any uniformly
bounded sequence of operators $\{ A_{k} \}$ on $\H_{n}$ which satisfies
$\lim_{k\rightarrow\infty}\| A_{k}\xi_{w} \| =0$ for all words $w$,
must converge in the strong operator topology to 0.
This completes the proof.
{\bx}

As a direct consequence of the lemma, one obtains the intuitive result that the powers of every non-unitary
contraction in $\gL_{n}$ converge in the weak operator topology to 0.
The lemma also yields deeper information on the structure theory of
contractions.

\begin{thm} The minimal isometric dilation of any non-unitary contraction in
$\gL_{n}$ is a shift.
\end{thm}

{\Pf}
In  \cite{nagfo} it is shown that every contraction has a
minimal isometric dilation. By the Wold decomposition for isometries, every
isometry is the orthogonal direct sum of a unitary and copies of the unilateral shift.
It is further shown in \cite{nagfo} that the powers of the adjoint of the
contraction
converging strongly to 0 (in other words, belonging to the class $C_{\cdot 0}$) is equivalent to
 the unitary part of its minimal
isometric dilation being vacuous. Hence the lemma yields the
result.
{\bx}

\begin{rem}
This appears to be new for $n=1$. The author could find no references, but it
is probably known in this case.

Recall that a contraction is called {\bf completely non-unitary} provided its
restriction to any non-zero reducing subspace is never unitary.
\end{rem}

\begin{cor} Every non-unitary contraction $L$ in $\gL_{n}$ is completely non-unitary. \end{cor}

{\Pf} Any non-zero reducing subspace for which the restriction of $L$ to it is unitary
would be contained in the unitary summand of the Wold decomposition for the minimal isometric dilation of $L$.
However, by the theorem this space is vacuous.
{\bx}

Another important result which comes out of the Sz.-Nagy and Foia\c{s}
machinery is that every completely non-unitary contraction possesses an
$\hinf$ functional calculus.

\begin{cor} \label{fnlcalc}
Every non-unitary contraction in $\gL_{n}$ has an $\hinf$ functional calculus.
\end{cor}
Given a non-unitary contraction $L$ in $\gL_{n}$, the collection of operators
defined by this $\hinf$ functional calculus is denoted  $\hinf(L)$.

The cardinality of the shift in the Wold decomposition of an isometry $V$ in
$\B (\H)$ is given by the dimension of $\H \ominus V\,\H$. For isometries in
the analytic Toeplitz algebra this cardinality can be both finite and infinite.
In fact if $\phi$ belongs to $\hinf$, the dimension of $H^{2}\ominus\phi\,H^{2}$ is finite exactly when the analytic inner symbol is continuous.
That is, when $\phi$ belongs to the disk algebra (see \cite{douglas}).
For the non-commutative algebras there is in general more room in the
orthogonal complement, and this cardinality turns out to always be infinite.
To prove this, first note the following result, the proof of which is
actually contained in the proof of Theorem 1.7 from \cite{topinv}.

\begin{thm} \label{inf_codim} For $n\geq 2$, if $A$ in $\gL_{n}$ has proper closed range, then
$\overline{{\rm Ran\,}(A)}$ has infinite codimension.
\end{thm}

\begin{rem} \label{pkcard}
For isometries $L$ in $\gL_{n}$,  the
infinite cardinality of $\H_{n}\ominus L\,\H_{n}$ is actually easy to see when
$(L\xi_{1},\xi_{1})=0$.
Indeed, let $P_{k}$ be the projection of $\H_{n}$ onto ${\rm span\,}\{\xi_{w}
: |w|=k\}$. Then,
\[ P_{k}\H_{n} \supseteq P_{k}({\rm Ran\,}(L)) = P_{k}\left(
\sum_{i=0}^{k-1}L(P_{i}\H_{n})\right). \]
The former space has dimension $n^{k}$, the latter has dimension at most
$\frac{n^{k}-1}{n-1}$. Summing over $k\geq 1$ proves the claim.
\end{rem}

In any event, it follows that the range of any non-outer operator has
infinite codimension. Recall that $A$ in $\gL_{n}$ is {\bf inner} if it is
an isometry and {\bf outer} if ${\rm Ran\,}(A)$ is dense in $\H_{n}$.

\begin{cor}
For $n\geq 2$, if $A$ in $\gL_{n}$ is not outer, then
$\overline{{\rm Ran\,}(A)}$ has infinite codimension.
\end{cor}

{\Pf}
Since $A$ is not outer, by the unique inner-outer factorization it can be written as $A=L\,B$ for some non-scalar isometry $L$
and outer operator $B$, both in $\gL_{n}$. But then,
\[ \overline{A\,\H_{n}} = \overline{L\,B\,\H_{n}} = L(\overline{B\,\H_{n}})
= L\,\H_{n}. \]
The latter space has infinite codimension by the theorem.
{\bx}

The general result can now be proved. The proof makes use of Fredholm theory
in $\gL_{n}$.

\begin{thm} Let $L$ be a non-unitary contraction in $\gL_{n}$, for $n\geq 2$. If  $V$ in $\B(\K)$ is its minimal
 isometric dilation, then
\[ \dim (\K\ominus V\,\K) = \infty. \]
\end{thm}

{\Pf}
The multiplicity of $V$ is given by the rank of $I-V\,V^{*}$.
From the construction of the minimal isometric dilation, this
is at least the rank of the operator $I-L\,L^{*}$. Suppose this number is
finite. Then $L^{*}$ has an essential left inverse, and hence $\ker L^{*} =
(\overline{{\rm Ran\,}L})^{\perp}$ is finite dimensional.
Thus by the previous corollary, $L$ must in fact be outer.

It now follows that the operators $L^{*}L$ and $L\,L^{*}$ are unitarily
equivalent. For, the partial isometry in the polar decomposition of $L$ is really
invertible and acts as the intertwining unitary. Therefore,
\[ {\rm rank\,}(I-L^{*}L) = {\rm rank\,}(I-L\,L^{*}) < \infty , \]
so that $L$ is an essential unitary.

As a Fredholm operator, $L$ has closed range and is thus surjective since it is outer. From \cite{topinv}, every operator in $\gL_{n}$ is injective. Hence $L$ is invertible. However,
it was also shown in \cite{topinv} that the essential norm of every operator in
$\gL_{n}$ is the same as its original norm. Since $\gL_{n}$ is inverse closed
\cite{topinv}, this implies that
\[ \|L^{-1}\| = \| L^{-1} \|_{e} = \| L^{*} \|_{e} = 1. \]
As an invertible isometry in $\gL_{n}$, $L$ must be scalar. This contradiction
completes the proof.
{\bx}

The investigation of the structure of $\H_{n}\ominus L\,\H_{n}$ will be
revisited  next section in the context of factoring.

In many ways contractions satisfying $(L\,\xi_{1},\xi_{1})=0$ are easier to
deal with. It is thus helpful to finish off this section by observing that
there is a large class of contractions in $\gL_{n}$ for which this inner
product is non-zero, however these operators are unitarily equivalent to
contractions in $\gL_{n}$ which have no scalar part. Observe that for any
operator $L$ in $\gL_{n}$, the space $\H_{n}\ominus \overline{L\,\H_{n}}$ belongs to
${\rm Lat\,}\gR_{n}^{*}$.

\begin{thm}
Suppose $L$ is a contraction in $\gL_{n}$ for which $\H_{n}\ominus \overline{L\,\H_{n}}$
contains an eigenvector of $\gR_{n}^{*}$. Then $L$ is unitarily equivalent
to a contraction in $\gL_{n}$ which has no scalar part. In addition, this
unitary implements an automorphism of $\gL_{n}$.
\end{thm}

{\Pf}
In \cite{topinv} the eigenvectors of $\gR_{n}^{*}$ are identified. Each
scalar $\lambda$ in the unit ball of $n$-dimensional Hilbert space defines
an eigenvector $v_{\lambda}$. Further, for each such vector there is a unitary
$U_{\lambda}$ on $\H_{n}$ for which ${\rm Ad\,}U_{\lambda}$ determines an
automorphism of $\gL_{n}$ with
\[ U_{\lambda}v_{\lambda} = \xi_{1}. \]

Thus, suppose some $v_{\lambda}$ belongs to $\H_{n}\ominus \overline{L\,\H_{n}}$. Then
the operator $U_{\lambda}L\,U_{\lambda}^{*}$ is a contraction in $\gL_{n}$ and
\[ (U_{\lambda}L\,U_{\lambda}^{*}\xi_{1},\xi_{1}) = (L\,v_{\lambda},v_{\lambda}) = 0, \]
which proves the result.
{\bx}

\section{Factoring}

In the analytic Toeplitz algebra $\gL_{1}=\gR_{1}=\T(\hinf)$, the associated
function theory yields a good factorization theory over the full algebra
(see \cite{douglas} and \cite{hoffman} for example). When moving to several
non-commutative variables, the strong link to the function theory is lost and
factorization becomes much more difficult to deal with. Nonetheless, positive
results such as  inner-outer factorization can be obtained.
Other factorization results can be obtained when norm restrictions are placed
on possible factors.

\begin{thm}
Let $w\in\F_{n}$. Over the unit ball of $\gL_{n}$, the isometry $L_{w}$ factors
only in the same way as the word $w$, modulo scalars in $\bbT$.
\end{thm}

{\Pf} Suppose $L_{w}=B\,C$ with $B$ and $C$ belonging to $b_{1}(\gL_{n})$.
It is clear that $C$ must in fact be an isometry.
For each $k\geq 1$, consider the corresponding form of (1) for  $B$ and $C$.
Let $\{b_{v}\}$ be the scalars and let $\{B_{v}\}$ be the operators for $B$ in this
decomposition. Use similar notation for $C$.

If $w=1$, then $I=b_{0}c_{0}I$. Whence, $|b_{0}|=|c_{0}|=1$ and the operators
$B$ and $C$ are scalar unitaries. Otherwise, put $w=z_{i_{1}}\cdots z_{i_{k}}$,
and note that $b_{0}c_{0}=0$.

First suppose that $b_{0}=0$. Then $B = \sum_{i=1}^{n}L_{i}B_{i}$ and equating factorizations
yields,
\[ L_{i_{1}}\left( \Bzis C\right) = L_{i_{1}}\left( L_{i_{2}}\cdots L_{i_{k}}\right) \]  and \[ L_{j}\left( B_{j}C\right) = 0 \,\,\,\, {\rm for\,}j\neq i_{1}. \]
In particular, $B_{j}C = 0$ for $j\neq i_1$. But every non-zero element of $\gL_{n}$ is injective \cite{topinv}, so that $B_{j}=0$
for $j\neq i_{1}$. Further, $B_{i_{1}}$ is a contraction and
$B_{i_{1}}C = L_{i_{2}}\cdots L_{i_{k}}$. Hence, by induction
one has $B_{i_{1}}=\lambda\,L_{u}$ and $C=\bar{\lambda}\,L_{v}$ where $uv=
z_{i_{2}}\cdots z_{i_{k}}$ and $\lambda$ belongs to $\bbT$. Thus, $B=\lambda\,
L_{i_{1}}L_{u}$ and $C=\bar{\lambda}\,L_{v}$.

Next suppose  $b_{0}\neq 0$, so that $c_{0}=0$.
Note first that
\begin{eqnarray*}
(\xi_{w},\xi_{z_{i_{k}}}) &=& (B\,C\xi_{1}, \xi_{z_{i_{k}}}) \\
                          &=& \sum_{v,\, u}b_{u}c_{v}\,(\xi_{uv},\xi_{z_{i_{k}}}) \\
                          &=&  b_{0}c_{i_{k}}+b_{i_{k}}c_{0} \\
                          &=& b_{0}c_{i_{k}}
\end{eqnarray*}
Hence if $w=z_{i_{k}}$, one would have $b_{0}c_{i_{k}} = 1$. As $B$ and $C$ are contractions
it would follow that $|b_{0}|=|c_{i_{k}}|=1$, and that $B$ is a scalar unitary.
Otherwise suppose $|w|>1$.
Inductively, one can show that
\[ 0=c_{0}=c_{i_{k}}=\ldots=c_{i_{2}\cdots i_{k}}. \]
To see this, observe that
$ 0 = (\xi_{w}, \xi_{z_{i_{k}}}) = b_{0}c_{i_{k}},$
 and hence $c_{i_{k}}=0$. Then
 suppose $0=c_{0}=c_{i_{k}}=\ldots=c_{i_{j}\cdots i_{k}}$ for some $j$,
$2<j\leq k$.
Equating Fourier coefficients of $L_{w}=B\,C$ shows that
\begin{eqnarray*}
0 &=& (\xi_{w}, \xi_{z_{j-1}\cdots z_{k}}) \\
  &=& (B\,C\,\xi_{1}, \xi_{z_{j-1}\cdots z_{k}}) \\
  &=&  b_{0}\czi + b_{i_{j-1}}c_{i_{j}\cdots i_{k}}+ \ldots +b_{i_{j-1}\cdots i_{k}}c_{0} \\
  &=& b_{0}\czi.
\end{eqnarray*}
Whence, $\czi=0$ as claimed. But then,
\begin{eqnarray*}
1 &=& (\xi_{w},\xi_{w}) \\
 &=& (B\,C\,\xi_{1},\xi_{w}) \\
  &=& b_{0}c_{w} + b_{i_{1}}c_{i_{2}\cdots i_{k}} + \ldots + b_{w}c_{0} \\
  &=& b_{0}c_{w}.
\end{eqnarray*}
Thus, $|b_{0}|=|c_{w}|=1$ and $B$ is a scalar unitary.
{\bx}

It is immediate that the generating isometries are irreducible over the
unit ball.

\begin{cor} For $n\geq 2$, each $L_{i}$ is irreducible over the unit ball of
$\gL_{n}$.
\end{cor}

In fact, many more isometries are irreducible over the unit ball of $\gL_{n}$.
Indeed, one can work harder to obtain the next result which includes a large collection of isometries.
For example, the isometries $L=\frac{1}{\sqrt{2}}(L_{1} + L_{2})$ and $L=\frac{1}{\sqrt{2}}L_{1}+\frac{1}{2}L_{2}^{2}+\frac{1}{2}L_{3}^{3}$ are irreducible.

\begin{thm} \label{irred1}
Suppose $L\sim\sum_{w\neq 1} a_{w}L_{w}$ is an isometry in $\gL_{n}$ for which
there is an $i$ with $a_{z_{i}}\neq 0$ and $R_{i} R_{i}^{*}(L\,\xi_{1})=a_{z_{i}}\xi_{z_{i}}.$ Then $L$ is irreducible over the unit ball of $\gL_{n}$, modulo
scalars in $\bbT$.
\end{thm}

{\Pf} Suppose $L = B\,C$ with $B$ and $C$ in $b_{1}(\gL_{n})$. The operator
$C$ must be  an isometry. As in the proof of the previous theorem,
consider the expansions of $B$ and $C$ determined by (1).
Recall that
\[ B = b_{0}I + \sum_{j=1}^{n}L_{j}B_{j} \,\,\,\,{\rm and}\,\,\,\,
   C = c_{0}I + \sum_{j=1}^{n}L_{j}C_{j}. \]
By equating unique factorizations  of $L = B\,C$ one obtains,
$b_{0}c_{0} = 0$. It is also clear that
\[ L_{i}\left( b_{0}C_{i} + \Bzi C \right) = L_{i}\left( a_{z_{i}}I
\right) . \]
First suppose $b_{0}=0$.
Then, $\Bzi C = a_{z_{i}}I \neq 0$.
Thus by the injectivity of all elements in $\gL_{n}$, the operator $B_{i}$ must be invertible.
Hence $C$ is also invertible. Therefore, as an
invertible isometry in $\gL_{n}$, $C$ must be a scalar unitary.

Next suppose $b_{0}\neq 0$ and so $c_{0}=0$. This corresponds to the case when
$B$ is a scalar unitary. Indeed, note first from the Fourier expansions,
\begin{eqnarray*} a_{z_{i}} &=& (L\,\xi_{1},\xi_{z_{i}}) \\
                           &=& (B\,C\,\xi_{1},\xi_{z_{i}}) \\
                            &=&  b_{0}c_{z_{i}}+ b_{z_{i}}c_{0} \\
                            &=&  b_{0}c_{z_{i}}.
\end{eqnarray*}
As $B$ is a contraction, $|c_{z_{i}}|\geq |a_{z_{i}}|$.
Now by hypotheses one has
\begin{eqnarray*}
(L-a_{z_{i}}L_{i})\,\xi_{1} &=& \sum_{j=1}^{n}R_{j}R_{j}^{*}(L-a_{z_{i}}L_{i})\,\xi_{1} \\
                            &=& \sum_{j\neq i}R_{j}R_{j}^{*}(L\,\xi_{1}).
\end{eqnarray*}
Further since $L=B\,C$,
\begin{eqnarray*}
B\,(C-c_{z_{i}}L_{i})\,\xi_{1} &=& \sum_{j=1}^{n} R_{j}R_{j}^{*}(L-c_{z_{i}}B\,L_{i})\,\xi_{1} \\
                               &=& \sum_{j\neq i}
R_{j}R_{j}^{*}(L\,\xi_{1}) + R_{i}R_{i}^{*}(L-c_{z_{i}}B\,L_{i})\,\xi_{1}.
\end{eqnarray*}
Evidently then,
\[ \|B\,(C-c_{z_{i}}L_{i})\,\xi_{1} \| \geq \| (L-a_{z_{i}}L_{i})\,\xi_{1}\|.\]
Thus the following is true:
\begin{eqnarray*}
1 = \| C\,\xi_{1} \| ^{2} &=& \| (C-c_{z_{i}}L_{i})\,\xi_{1} \|^{2} + |c_{z_{i}}|^{2} \\
                          &\geq& \| B\,(C-c_{z_{i}}L_{i})\,\xi_{1}\|^{2} + |c_{z_{i}}|^{2} \\
                          &\geq& \| (L-a_{z_{i}}L_{i})\,\xi_{1}\|^{2} + |a_{z_{i}}|^{2} \\
                          &=& \| L\,\xi_{1} \| ^{2} = 1.
\end{eqnarray*}
Therefore, $|a_{z_{i}}|=|c_{z_{i}}|$ and $|b_{0}|=1$, which shows that  $B$ is a scalar
unitary.
{\bx}

\begin{rem}\label{irred2}
This proof can be perturbed to include broader classes of isometries. For example, any isometry $L$ which satisfies,
\[ R_{i}R_{i}^{*}(L\,\xi_{1}) = a_{wz_{i}}\xi_{wz_{i}}, \]
for some $i$ and word $w$ is irreducible over the unit ball of $\gL_n$.
\end{rem}

As it turns out, the unique factorizations over the unit ball of $\gL_{n}$
discussed above do not hold over the full algebra. Remarkably, even the operator $L_{2}$
has  proper factorizations in $\gL_{n}$. This comes out of an interesting
result from the function theory.

\begin{lem} \label{fntheory} The function $f\,(z) = \sum_{k\geq 0}\frac{z^{k}}{k+1}$
belongs to $H^{2}\setminus\hinf$. However, $\onef$ defines a function which lies in
$\hinf$.
\end{lem}

{\Pf}
Since $f$ is analytic on the unit disk and the Fourier coefficients of $f$ are $\ell_{2}$-summable,
 the function belongs to $H^{2}$ with
\[ \| f \|_{2}=\left(\sum_{k\geq 0}\frac{1}{(k+1)^2}\right)^{
\frac{1}{2}}. \]
For $|z|<1$, $f\,(z)$ is defined by the formula
\begin{eqnarray} \label{log}
z\, f\,(z) = -\log (1-z)
\end{eqnarray}
for the principal branch of the logarithm. Given $r$ such that $0<r<1$, let
$f_{r}$ be the function on $\bbT$ defined by $f_{r}(e^{i\theta})=f\,(r\,e^{i
\theta})$. Then the identity
\[ \lim_{r\rightarrow 1^{-}}f\,(r) = \lim_{r\rightarrow 1^{-}}\frac{-\log(1-r)}{r} = \infty, \]
together with the continuity of $f$ on the disk, shows that
\[ \| f \|_{\hinf} := \lim_{r\rightarrow 1}\|f_{r}\|_{\infty} = \infty. \]
Hence the function $f$ is not in $\hinf$.

To prove that $\onef$ defines a function in $\hinf$,  it is required to show that $\onef$
defines an analytic function on the unit disk with $\|\onef\|_{\infty}<\infty$. Now,  $\onef$ is analytic
on the unit disk by (2).
To see that $f$ is bounded below first observe the
identity $1-e^{i\,\theta} = (2 \sin \frac{\theta}{2})\,e^{
\frac{\theta - \pi}{2}i}$.
 Hence for $|\theta|\leq\pi$ with $\theta\neq 0$,
\begin{eqnarray*}
|f\,(e^{i\,\theta})|^{2} &=& \left| \frac{-\log (1-e^{i\,\theta})}{e^{i\,
\theta}} \right|^{2} \\
         &=& \left| \frac{\log |2 \sin \frac{\theta}{2}| + i\,\left(\frac{\theta
-\pi}{2}\right)}{e^{i\,\theta}}\right|^{2} \\
         &=& \left(\log |2 \sin \frac{\theta}{2} | \right)^{2} + \left( \frac{
\theta-\pi}{2}\right)^{2}.
\end{eqnarray*}
But $\left(\frac{\theta-\pi}{2}\right)^{2}\geq \frac{\pi^{2}}{16}$ for
$-\pi\leq\theta\leq\frac{\pi}{2}$ and $|2 \sin \frac{\theta}{2}|\geq\sqrt{2}$
for $\frac{\pi}{2}\leq\theta\leq\pi$. Thus,
\[ |f\,(e^{i\,\theta})|^{2} \geq \min \left\{ \frac{(\log 2)^{2}}{4}, \frac{
\pi^{2}}{16} \right\}  = \frac{(\log 2)^{2}}{4} \]
for $\theta\neq 0$. It now follows that
$\|\onef\|_{\infty}<\infty$, and the proof is finished.
{\bx}

This unusual function theoretic result allows one to construct explicit
factorizations which are exclusive to the non-commutative setting.

\begin{thm} Suppose $L$ is an isometry in $\gL_{n}$ for which $\H_{n}\ominus
L\,\H_{n}$ contains the range of an isometry $X$ in $\gL_{n}$. Then $L$ has
proper factorizations in $\gL_{n}$.
\end{thm}

{\Pf} As the ranges of the isometries $X^{k}L$ are pairwise orthogonal for
$k\geq 0$, the operator $A=\sum_{k\geq 0}\frac{1}{k+1}X^{k}L$ belongs to
$\gL_{n}$ with
\[ \| A  \| = \left( \sum_{k\geq 0}\frac{1}{(k+1)^{2}} \right) ^{\frac{1}{2}}.\]
Let $g= \frac{1}{f}$ be the $\hinf$ function obtained in the previous lemma. By the
$\hinf$ functional calculus for $X$ (Corollary ~\ref{fnlcalc}), $g\,(X)$
defines an operator in $\gL_{n}$. The claim is that $g\,(X)\,A=L$.
This comes as a result of a more general fact.

Note that given $h$ in $\htwo$, an operator $h\,(X)\,L$ can be naturally defined in
$\gL_{n}$. Indeed, one can set
\[ h\,(X)\,L = \sum_{k\geq 0} \hht \,(k)\,X^{k}\,L. \]
where the $\hht(k)$ are the  Fourier coefficients for $h$.
Clearly, the map from $\htwo$ to $\B(\H_{n})$ which sends $h$ to $h\,(X)\,L$
is isometric. The key is that this map is also continuous from the topology of
weak vector convergence in $\htwo$ to the $\WOT$ in $\gL_{n}$. To see this,
suppose $h_{m}$ converges weakly to $h$ in $\htwo$. Without loss of generality
assume $h=0$. Then $\hht_{m}(k)$ converges to 0 for each $k$ and
\[ \sup_{m}\|h_{m}\|_{2}=c<\infty , \]
for some constant $c$.
Let $x$ and $y$ be unit vectors in $\H_{n}$ and let $S_k$ be the orthogonal projection onto ${\rm Ran}\,(X^k L)$. Then
\begin{eqnarray*}
\sum_{k\geq 0} |(X^k Lx,y)|^2 &=& \sum_{k\geq 0} |(X^k Lx,S_k y)|^2 \\
                              &\leq& \sum_{k\geq 0} ||S_k y||^2 \\
                              &\leq& ||y||^2 = 1.
\end{eqnarray*}
Thus, given $\varepsilon >0$ one can choose $N=N(\varepsilon)$ for which the $N$th $\ell_{2}$ tail of the above series is smaller than $\varepsilon$.
Then for each $m$ the Cauchy-Schwarz inequality shows that
\begin{eqnarray*}
|(h_{m}(X)\,L\,x,y)| &=& \left| \sum_{k\geq 0}\hht_{m}(k)\,(X^{k}Lx,y) \right| \\
               &\leq& \sum_{0\leq k \leq N}|\hht_{m}(k)| +
\left(\sum_{k>N}|\hht_{m}(k)|^{2}\right)^{\frac{1}{2}} \varepsilon \\
               &\leq& \sum_{0\leq k \leq N}|\hht_{m}(k)| + c\,\varepsilon .
\end{eqnarray*}
As $\hht_{m}(k)$ converges to 0 for each $k$, it follows that $h_{m}(X)\,L$
converges $\WOT$ to 0.

Recall that the analytic trigonometric polynomials are weak* dense in $\hinf$
\cite{douglas}. Let ${g_{m}}$ be such a sequence converging weak* to $g$.
From the definition of this weak* topology, the sequence ${g_{m}}$ converges
weakly to $g$ in $\hinf$. Thus, the sequence ${g_{m}f}$ converges weakly to
$g\,f=1$ in $\htwo$. Hence since each $g_{m}$ is a polynomial,
\[ g_{m}(X)\,A = g_m(X)\left( f(X)L\right) = (g_{m}f)(X)\,L \stackrel{\WOT}{\longrightarrow} (g\,f)(X)\,L =L.
\]
Further, by the $\hinf$ functional calculus for $X$, $g_{m}(X)$ converges
$\WOT$ to $g\,(X)$ \cite{nagfo}. Therefore,
\[ g_{m}(X)\,A \stackrel{\WOT}{\longrightarrow}
 g\,(X)\,A. \]  Whence $g\,(X)\,A=L$.

It remains to observe that $g\,(X)$ and $A$ are both not invertible. The
invertibility of $g\,(X)$ in $\gL_{n}$ would imply the invertibility of $g$ in
$\hinf$, contradicting the previous lemma. If $A$ was invertible, it would
be the scalar multiple of an invertible isometry in $\gL_{n}$, hence scalar
itself by Proposition ~\ref{unit_prop}. The proof is now complete.
{\bx}

\begin{rem}
The theorem really is exclusive to the non-commutative setting. The hypothesis of the theorem cannot be satisfied when $n=1$.  For if $\phi$ and
$\psi$ are inner functions in $\hinf$, then the function $\phi\,\psi=\psi\,
\phi$ belongs to $\phi\,H^{2}$ and $\psi\,H^{2}$.
\end{rem}

As a surprising consequence of the theorem the reducibility of the generating isometries
is revealed in the non-commutative setting.

\begin{cor} For $n\geq 2$, each $L_{i}$ has proper factorizations over
$\gL_{n}$.
\end{cor}

{\Pf} The isometry $X$ can be taken to be $L_{j}$ for $j\neq i$.
{\bx}

In fact, there is a large collection of isometries which can be seen to be
reducible in this manner. Note that by unique factorization, every operator
$L$ in $\gL_{n}$ with $(L\,\xi_{1},\xi_{1})=0$ can be written as
\[ L = \sum_{|w|=k}L_{w}A_{w}, \]
for some $k\geq 1$ and $A_{w}$ in $\gL_{n}$.

\begin{cor} Let $L$ be an isometry in $\gL_{n}$ with $L=\sum_{|w|=k}L_{w}A_{w}$, for some $k\geq 1$ and $A_{w}$ in $\gL_{n}$. Suppose either one of the following conditions holds:
\begin{enumerate}
\item  There is an  $A_{w}$ for which there exists a $B\neq 0$ in $\gL_{n}$
       with the range of $B$ orthogonal to that of $A_{w}$, or
\item  Some $A_{w}$ is a scalar multiple of an isometry.

\end{enumerate}
Then $L$ properly factors in $\gL_{n}$.
\end{cor}

{\Pf} Let $\xi$ and $\eta$ be vectors belonging to $\H_{n}$ throughout the proof.
To prove (1), set $A=L_{w}B$.  Write the inner-outer factorization for $A$ as
$A=L_{A}C$. Since $C$ is outer, there are vectors $\zeta_{m}$ such that
$\eta=\lim_{m\rightarrow\infty}C\,\zeta_{m}$. It follows that
\begin{eqnarray*}
(L\,\xi,L_{A}\eta) &=& \lim_{m\rightarrow\infty}(L\,\xi,L_{A}C\,\zeta_{m}) \\
              &=& \lim_{m\rightarrow\infty}(L\,\xi,L_{w}
                  B\,\zeta_{m}) \\
              &=& \lim_{m\rightarrow\infty} (L_{w}A_{w}\xi,L_{w}B\zeta_{m}) \\
              &=& \lim_{m\rightarrow\infty} (A_{w}\xi,B\zeta_{m}) \\
              &=& 0.
\end{eqnarray*}
Thus, $L_{A}$ is an isometry with range orthogonal to the range of $L$.

Lastly, suppose that $A_{w}$ is a scalar multiple of an
isometry. Write $L$ as $L=L_{w}A_{w} + A$. Let $A=L_{A}B$ be the inner-outer
factorization of $A$. Since
\[ \overline{A\,\H_{n}} = \overline{L_{A}B\,\H_{n}} = L_{A}\H_{n}, \]
the ranges of $L_{w}$ and $L_{A}$ are orthogonal.
Further,  one has
\[ A_{w}^{*}A_{w} + B^{*}B = I. \]
So $B$ is also a scalar multiple of an isometry, and it is therefore scalar
since it is outer. Suppose that $A_{w}^{*}A_{w}=a^{2}I$ and $B=\beta\,I$. Let $c = a^{2}|\beta|^{-2}$, and let $X$ be the operator
\[ X=L_{w}A_{w} - c\,\beta\,L_{A}. \]
Then $X$ is an operator with range orthogonal to the range of $L$. Indeed,
\begin{eqnarray*}
(L\,\xi,X\,\eta) &=& (L_{w}A_{w}\,\xi,L_{w}A_{w}\,\eta)-(A\,\xi,c\,\beta\,L_{A}\,\eta) \\
            &=& a^{2}(\xi,\eta)-c\,|\beta|^{2}(\xi,\eta) \\
            &=& 0.
\end{eqnarray*}
 As in the proof of (1), using the
inner part of $X$ yields a desired isometry.
{\bx}

It is worthwhile to point out a striking special case of the first condition
in the corollary.

\begin{cor} Let $L$ be an isometry in $\gL_{n}$ with $L=\sum_{|w|=k}L_{w}A_{w}$, for some $k\geq 1$ and $A_{w}$ in $\gL_{n}$. If any $A_{w}=0$, then $L$
properly factors in $\gL_{n}$.
\end{cor}

\begin{rem} Obviously there are many isometries which satisfy
the first condition. Further, all of the isometries shown to be irreducible over the unit ball of the algebra in Theorem ~\ref{irred1} and Remark ~\ref{irred2} are reducible over the full algebra since they satisfy the second condition.
Other isometries which satisfy the second condition include the collection of
all operators which are
the sum of pairwise orthogonal words. For in this case every non-zero
$A_{w}$ would necessarily be a scalar multiple of an isometry.

There are also other more specialized classes of isometries which can be
factored using this method. As an example, let $f$ and $g$ belong to $\hinf$
with \[ |f|^{2} + |g|^{2} = 1 \] on $\bbT$. Such functions can be found by using the
logmodularity of $\hinf$ \cite{hoffman}.
Then
\[ L = L_{1}\,f\,(L_{1}) + L_{2}\,g\,(L_{1}) \]
is an isometry in $\gL_{2}$ which apparently does not  satisfy  the
conditions in the corollary. Let $\alpha = f\,(0)$ and $\beta = g\,(0)$ and choose $\lambda$ in $\bbT$ such that
\[ \lambda( \alpha\,\overline{\beta}) = \overline{\alpha}\,\beta.
\]  Then $L$ and the isometry
\[ X = \frac{1}{|\alpha|^{2}+|\beta|^{2}} \left( \beta\,L_{1}L_{2} -
\lambda\,\alpha\,L_{2}^{2} \right) \]
have orthogonal ranges. Indeed, for $\xi$ and $\eta$ in $\H_{n}$ one has
\begin{eqnarray*}
(|\alpha|^{2} + |\beta|^{2})\left( L\,\xi,X\,\eta \right) &=&
\left( f\,(L_{1})\,\xi,\beta\,L_{2}\,\eta\right) - \left( g\,(L_{1})\,\xi,\lambda\,
\alpha\,L_{2}\,\eta\right) \\
&=& \alpha\,\overline{\beta}\,\left( \xi,L_{2}\,\eta\right) - \beta\,\overline{\alpha}\,
\overline{\lambda}\,\left( \xi,L_{2}\,\eta \right) \\
&=& 0.
\end{eqnarray*}

\end{rem}

The reducibility of this large collection of isometries, together with the fact that
the orthogonal complement of the range is always infinite dimensional (Theorem  ~\ref{inf_codim}), leads
one to believe that perhaps the theorem can be applied to every isometry $L$
with $(L\,\xi_{1},\xi_{1})=0$. However, this is not the case. The trouble is
that difficulties arise when the space $\H_{n}\ominus L\,\H_{n}$ is too `thin'
at each level of the $\H_{n}$ tree. That is, the dimension of
$P_{k}(\H_{n}\ominus L\,\H_{n})$ remains small as $k$ increases (recall
Remark ~\ref{pkcard}).

\begin{thm} There are isometries $L$ in $\gL_{n}$ with $(L\,\xi_{1},\xi_{1})=0$ for which $\H_{n}\ominus L\,\H_{n}$ does not contain the range of an isometry
in $\gL_{n}$.
\end{thm}

{\Pf} For $k\geq 0$, put
\[ x_{k} = R_{1}^{k}R_{2} \sum_{|w|=k} R_{w}\xi_{w}, \]
and let $x$ be the unit vector
\[ x = \sum_{k\geq 0} 2^{\frac{-k-1}{2}}\| x_{k} \|^{-1}\, x_{k}. \]
Suppose $y$ is in $\H_{n}$ with \[ (R_{u}^{*}x,y)=0 \] for all words $u$ in $\F_{n}$.
Now given $k\geq 0$, choose a word $u$ with $|u|=k$. Then
\begin{eqnarray*}
0 &=& \left( 2^{\frac{k+1}{2}}\| x_{k} \|^{1} \, R_{u}^{*}R_{2}^{*}(R_{1}^{k})^{*} x, y \right) \\
  &=& ( R_{u}^{*}\bigl( \sum_{|w|=k}R_{w} \xi_{w}\bigr), y ) \\
  &=& \left( \xi_{u},y \right).
\end{eqnarray*}
Therefore, $y=0$.

Next, write $x$ as $x = L\,\eta$, where $L$ is an isometry in
$\gL_{n}$ and $\eta$ is an $\gR_{n}$ cyclic vector (every vector in $\H_{n}$
can be written in this form \cite{topinv}). The claim is that $L$ is the desired isometry.
As $(\eta, \xi_{1})\neq 0$ and $(x,\xi_{1})=0$, one has
\mbox{$(L\,\xi_{1}, \xi_{1})=0$}. Suppose $X$ is an isometry in $\gL_{n}$ with
range contained in $\H_{n}\ominus L\,\H_{n}$. Then the vectors
$X\,\xi_{u} = X\,(R_{u}\,\xi_{1}) = R_{u}\,( X\,\xi_{1})$ are orthogonal to $L\,\eta = x$ for every $u$ in $\F_{n}$.
In other words,
\[ \left( X\,\xi_{1},R_{u}^{*}x \right) = 0 \]
for every word $u$. Thus, by the above argument one would have $X\,\xi_{1}=0$, whence $X=0$. This contradiction completes the proof.
{\bx}

So this method cannot be applied to all isometries in $\gL_{n}$.
Nevertheless, with the large body of examples it is still reasonable to
make the guess that every isometry $L$ in $\gL_{n}$ with $(L\,\xi_{1},
\xi_{1})=0$ properly factors over the full algebra.

\section{Ideals and Invariant Subspaces}

The characterization of the $\WOT$-closed right and two sided ideals (${\rm Id}
_{r}(\gL_{n})$ and ${\rm Id}\,(\gL_{n})$) in
\cite{topinv} and \cite{topalg} is complete.
The main theorem from \cite{topalg} is stated as follows.

\begin{thm} Let $\mu\,:\,{\rm Id}_{r}(\gL_{n})\,\rightarrow\,\Lat (\gR_{n})$
be given by $\mu\,(\I)= \overline{\I\,\xi_{1}}$. Then $\mu$ is a complete
lattice isomorphism. The restriction of $\mu$ to the set ${\rm Id}\,(\gL_{n})$
is a complete lattice isomorphism onto $\Lat (\gL_{n}) \cap \Lat (\gR_{n})$.
The inverse map $i$ sends a subspace $\M$ to
\[ i\,(\M) = \{ J\,\in\,\gL_{n}\,:\,J\,\xi_{1}\,\in\,\M\}. \]
\end{thm}

The maps $\mu$ and $i$ are still well defined when considering
${\rm Id}_{l}(\gL_{n})$ and $\Lat (\gL_{n})$.
Although technical difficulties were encountered by the authors, but a similar
characterization was expected for left ideals.
The key observation for right and two-sided ideals is that the subspace $\mu\,(
\I)$ is the
full range of the ideal $\I$. Indeed, $\overline{\I\,\xi_{1}} = \overline{\I\,
\gL_{n}\xi_{1}}= \overline{\I\H_{n}}. $
This is not true for left ideals, and is why the methods of the authors cannot
be applied in this setting.

Towards the identification of right ideals  it is first proved that $\mu\,i =
{\rm id}$. It is then shown that this leads to the conclusion,
$\overline{\I\, \xi} = \overline{i\,\mu\,(\I)\, \xi}$ for every vector $\xi$ in
$\H_{n}$.
The proof that $i\,\mu = {\rm id}$
exploits this fact together with
 the following more general result about ideals in $\gL_{n}$.

\begin{prop}
Let $\I_{1}$ and $\I_{2}$ both be $\WOT$-closed right, left or two-sided ideals
in $\gL_{n}$. If $\overline{\I_{1}\xi}=\overline{\I_{2}\xi}$ for all $\xi$ in
$\H_{n}$
, then $\I_{1}=\I_{2}$.
\end{prop}

{\Pf} In \cite{topinv} it was shown that the weak* and weak operator
topologies on $\gL_{n}$ coincide.
Suppose that $\phi$ is a $\WOT$-continuous functional on $\gL_{n}$
which annihilates $\I_{1}$. Then again from \cite{topinv}, there are vectors $\xi$ and
$\eta$ in $\H_{n}$ such that
\[ \phi\,(J) = \left( J\, \xi , \eta \right) \]
for $J$ in $\gL_{n}$. But then $\eta$ is orthogonal to $\overline{\I_{1}\xi}
= \overline{\I_{2}\xi}$, and hence $\phi$ annihilates $\I_{2}$ as well.
Repeating the argument by exchanging the r\^{o}les of $\I_{1}$ and $\I_{2}$
shows that the two ideals are identical.
{\bx}

\begin{rem}
Now, let $\I_{1}$ and $\I_{2}$ be $\WOT$-closed left ideals of $\gL_{n}$.
Notice that, $\overline{\I_{1}\xi_{1}} = \overline{\I_{2}\xi_{1}}$
implies $\overline{\I_{1}\xi} = \overline{\I_{2}\xi}$ when $\xi=R\,\xi_{1}$ for
some
 isometry $R$ in $\gR_{n}$.
 Indeed, one would have
\[ \overline{\I_{1}\xi} = \overline{\I_{1}R\,\xi_{1}}=R\,\overline{\I_{1}\xi_{1}
}
= R\,\overline{\I_{2}\xi_{1}} = \overline{\I_{2}R\,\xi_{1}} = \overline{\I_{2}
\xi}.  \]
These vectors form a dense collection of vectors in $\H_{n}$. Whether this
implies the same is true for all vectors in $\H_{n}$ is unclear.
As Remark ~\ref{unbded}  points out, this requires an understanding of unbounded
 $\WOT$-convergence. Thus, it becomes apparent that there are
difficulties encountered when considering left ideals.
\end{rem}

Upon further investigation concrete differences become evident. In particular,
the $\WOT$-closed
right ideal generated by a finite collection of isometries with pairwise
orthogonal ranges is exactly the algebraic right ideal they generate. For
left and two-sided ideals the corresponding result turns out to be false even   for one isometry with norm
closure.

\begin{thm}
The algebraic two-sided ideal in $\gL_{n}$ generated by $L_{2}$ is not norm closed.
\end{thm}

{\Pf} The operator \[ A=\sum_{k\geq 0}\frac{1}{k+1}L_{1}^{k}L_{2} \] clearly belongs to the norm closure of the algebraic two-sided (in fact left) ideal generated by $L_{2}$.  Suppose that $A$ could be  written as
\[ A=  \sum_{i=1}^{p}B_{i}L_{2}C_{i} \]
with each $B_{i}$ and $C_{i}$ in $\gL_{n}$. Put $B_{i}\sim\sum_{w}b_{w}^{i}L_{w}$ and $C_{i}\sim\sum_{w}c_{w}^{i}L_{w}$. Then the unique factorization
in $\gL_{n}$ shows for each $k\geq 0$,
\begin{eqnarray*}
 \frac{1}{k+1} &=& (A\,\xi_{1}, \xi_{z_{1}^{k}z_{2}}) \\
               &=& \sum_{i=1}^{p}(B_{i}L_{2}C_{i}\xi_{1},\xi_{z_{1}^{k}z_{2}}) \\
               &=& \sum_{i=1}^{p}b_{z_{1}^{k}}^{\,i}\,c_{0}^{i}.
\end{eqnarray*}
By compressing the operators $B_{i}$ to the subspace
${\rm span}\,\{\xi_{z_{1}^{k}}:k\geq 0\}$,  one sees that each of the operators \[ h_{i}(L_{1}) = \sum_{k\geq 0}b_{z_{1}^{k}}^{\,i}L_{1}^{k} \,\,\,\,{\rm for}\,
1\leq i \leq p \] must be in  $\hinf(L_{1})\stackrel{\sim}{=}\hinf$.
Hence the function \[ \sum_{k\geq 0}\frac{z^{k}}{k+1} = \sum_{i=1}^{p}c_{0}^{i}h_{i} \]
would belong to  $\hinf$,
a contradiction (see Lemma ~\ref{fntheory}). Therefore $A$ does not belong to the algebraic two-sided ideal
generated by $L_{2}$.
{\bx}

As an immediate corollary of the proof, the corresponding fact about left ideals is proved.

\begin{cor}
The algebraic left ideal $\gL_{n}L_{2}$ is not norm closed.
\end{cor}

{\Pf} Consider the same operator $A$. Simply use the proof of the theorem with
$p=1$ and $C_{1}=I$.
{\bx}

\begin{rem}
It seems reasonable to expect that the norm and $\WOT$ closures of the preceding ideals are distinct. It also becomes apparent that proving this
would be quite subtle. Indeed, it is difficult to construct a bounded
operator belonging to
 the weak closure of $\gL_{n}L_{2}$ without being in the norm closure of
$\gL_{n}L_{2}$.
\end{rem}

Even with these differences it is still surprising that the analogous
identification of left ideals does not hold.
It turns out that  the subspaces $\M$ in $\Lat (\gL_{n})$ for which $\mu\,i\,(\M
) =
\M$ do not fill out the entire subspace lattice.

\begin{thm}
There exists $\M\neq\{ 0 \}$ in $\Lat (\gL_{n})$ for which the associated left
ideal $i\,(\M)$ is trivial.
\end{thm}

{\Pf} Define an isometry $R$ in $\gR_{n}$ by
\[ R = \sum_{k\geq 0}\lambda_{k}\, R_{1}^{k}R_{2}, \]
where the scalars $\lambda_{k}$ satisfy $\sum_{k\geq 0}|\lambda_{k}|^{2}=1$ but
$\sum_{k\geq 0}\lambda_{k}z^{k}$ is not in $\hinf$. For example,
\[ \lambda_{k} = \frac{c}{k+1} \,\,\,\,{\rm where}\,\,\,\, c=\left(\sum_{k\geq 1
}\frac{1}{k^{2}}\right)^{- \frac{1}{2}}. \]
 Let $\M$ be the subspace in $\Lat (\gL_{n})$
given by $\M=R\,\H_{n}$. Actually, every cyclic $\gL_{n}$-invariant subspace is
of this form for some isometry in $\gR_{n}$ (see \cite{po}, \cite{topinv},
\cite{po1} and \cite{po2}). Now,
\[ i\,(\M)= \{ J\in\gL_{n}:J\,\xi_{1}\in R\, \H_{n}\}.\] Suppose there is a
non-zero $J$ in $\gL_{n}$ and $\xi$ in $\H_{n}$ for which $J\,\xi_{1} = R\,\xi$.
Put $\xi
=\sum_{w}a_{w}\xi_{w}$ and let $v$ be a word of minimal length such that
$a_{v}\neq 0$. Since $R\,\xi = \sum_{w}a_{w}L_{w}(R\,\xi_{1})$ one has
\[ J \sim \sum_{|w|\geq|v|,\, k\geq 0} a_w\lambda_k\,L_wL_2L_1^k
 . \]
Let $Q$ be the projection onto the subspace ${\rm span}\,\{\xi_{z_{1}^{k}}:k\geq
 0\}$.
Evidently then,
\[ Q\,L_{2}^{*}L_{v}^{*}J\,Q = a_{v}\,\sum_{k\geq 0}\lambda_{k}\,L_{1}^{k}Q, \]
incorrectly implying that $J$ would be unbounded.
Therefore, it follows that $i\,(\M)=0$.
{\bx}

There is still a strong relation between ${\rm Id}_{l}(\gL_{n})$ and
$\Lat (\gL_{n})$. Essentially, it is determined by those isometries in
$\gR_{n}$ which do not have the qualities of those used in the proof of the
theorem.

\begin{defn}
An isometry $R$ in $\gR_{n}$ is called a {\bf flip} if there is a non-zero
vector $\xi$ in $\H_{n}$ and an operator $J$ in $\gL_{n}$ with
\[ R\,\xi = J\,\xi_{1}. \]
Call $R$ a {\bf cyclic flip} if there are operators $J_{\alpha}$ in $\gL_{n}$
such that $J_{\alpha}\xi_{1}\,\in\,R\,\H_{n}$ and
\[ R\,\xi_{1} = \lim_{\alpha}J_{\alpha}\xi_{1}. \]
\end{defn}

The motivation for these definitions is when $R\,\xi_{1}= J\,\xi_{1}$,
for some $J$ in $\gL_{n}$. This means that the Fourier coefficients of $R$ can
be `flipped' into an element of $\gL_{n}$.

\begin{prop}
Let $R$ be an isometry in $\gR_{n}$ with $\M=R\,\H_{n}$. The following are
equivalent:
\begin{enumerate}
\item $R$ is a flip,
\item $i\,(\M)\neq 0$.
\end{enumerate}
\end{prop}

{\Pf} This is straight from the definitions of $i\,(\M)$ and flip isometries.
{\bx}

\begin{rem}
It should be noted that when $R$ is a flip, the subspace
 $\mu\,i\,(\M)$ is `large'. Indeed, suppose
$J\neq 0$ belongs to $i\,(\M)$. Then there is a vector
$\xi$ in $\H_{n}$ with $R\,\xi=J\,\xi_{1}$.  Write $\xi$ as
$\xi=S\,\eta$, where $S$ is an isometry in $\gR_{n}$ and $\eta$ is an $\gL_{n}$-
cyclic vector (this can be done for any vector in $\H_{n}$ \cite{topinv}).
Note that the set $\gL_{n}J$ is contained in $i\,(\M)$. Thus,
\[ R\,S\,\H_{n}= \overline{\gL_{n}R\,S\,\eta} = \overline{\gL_{n}R\,\xi}
= \overline{\gL_{n}J\,\xi_{1}} \subseteq \mu\,i\,(\M), \]
which shows that $\mu\,i\,(\M)$ contains the range of the isometry $R\,S$.
\end{rem}

\begin{prop} \label{lattice}
Let $R$ be an isometry in $\gR_{n}$ with $\M=R\,\H_{n}$.
The following are equivalent:
\begin{enumerate}
\item $R$ is a cyclic flip,
\item $\mu\,i\,(\M) = \M$.
\end{enumerate}
\end{prop}

{\Pf} It is always true that $\mu\,i\,(\M)\subseteq\M$.
Suppose that $R$ is a cyclic flip, so that $R\,\xi_{1}=\lim_{\alpha}J_{\alpha}\xi_{1}$ with each
$J_{\alpha}\xi_{1}\,\in\,\M\cap\gL_{n}\xi_{1}$.
Then $R\,\xi_{1}$ lies in $\mu\,i(\M)$, and hence
\[ \M = R\,\H_{n} = \overline{\gL_{n}R\,\xi_{1}} \subseteq \mu\,i\,(\M). \]
If this latter inclusion holds, then $R\,\xi_{1}$ is such a limit since
\[ R\,\xi_{1}\in \M = \mu\,i\,(\M) = \overline{\{J\in\gL_{n}:J\,\xi_{1}\in R\,\H
_{n}\}\,\xi_{1}}.\]
{\bx}

This gives a good one vector characterization of cyclic subspaces $\M$ in
$\Lat (\gL_{n})$ for which $\mu\,i\,(\M)=\M$.
The following corollary shows that the above condition is satisfied for a
wealth of examples. For instance, consider the situation below even for
$\xi=\xi_{1}$.

\begin{cor} \label{lat_cor}
Let $\M=R\,\H_{n}$, where $R$ is an isometry in $\gR_{n}$. If there is an
$\gL_{n}$-cyclic vector $\xi$ such that $R\,\xi=J\,\xi_{1}$ for some $J$ in
$\gL_{n}$, then $\mu\, i\,(\M)=\M$. Further, if $J=L_{w}$ for some word $w$,
then $i\,(\M)$ is exactly the $\WOT$-closed left ideal generated by $L_{w}$.
\end{cor}

{\Pf} The condition in the previous proposition is satisfied since,
\[ R\,\xi_{1}\in R\,\H_{n} = \overline{\gL_{n}R\,\xi} = \overline{\gL_{n}J\,\xi_{1
}}, \]
which shows that $R$ is a cyclic flip.

In general, given $A$ in $\gL_{n}$,
\[ A\,J\,\xi_{1} = A\,R\, \xi = R\,A\, \xi\,\in\,R\,\H_{n} = \M. \]
Thus, the $\WOT$-closed left ideal generated by $J$ is contained in $i\,(\M)$.
The other inclusion is true for the words $J=L_{w}$. Indeed, in this case it is easy to see that
every $A$ in $i\,(\M)$ has a Fourier expansion of the form
\[ A \sim \sum_{v\in\F_{n}}a_{vw}L_{vw} . \]
Actually, an analogous claim can be made for the right and two-sided
$\WOT$-closed ideals generated by $L_{w}$.
However, in \cite{topinv} it was shown that the Cesaro sums for $A$,
\[ \Sigma_{k}(A) = \sum_{|v|<k}\left( 1-\frac{|v|}{k} \right)a_{v}L_{v}\,\in\,\gL_{n}L_{w} \]
converge in the strong* topology to $A$. Hence $i\,(\M)$ is contained in the
$\WOT$-closed left ideal generated by $L_{w}$.
{\bx}

It has been mentioned that every vector $J\,\xi_{1}$ with $J$ in $\gL_{n}$
factors as $J\,\xi_{1}=R\,\xi$ for some isometry $R$ in $\gR_{n}$ and
$\gL_{n}$-cyclic vector $\xi$.
As the corollary observes,  one has that $i\,(R\,\H_{n})$ always contains the
$\WOT$-closed left ideal generated by $J$. The other inclusion holds for words $J=L_{w}$. Proving the other inclusion holds
in full generality would require an understanding of unbounded $\WOT$-convergence in these algebras. This is discussed further below.

\begin{rem} \label{unbded}
The corollary shows that the image  in $\Id_{l}(\gL_{n})$ contains
the left ideals generated by the words $L_{w}$.
It is not clear whether the image in $\Id_{l}(\gL_{n})$ is surjective.
Given a $\WOT$-closed left ideal $\I$, it is always true that $\I\subseteq\,i\,
\mu\,(\I)$. In general, the other inclusion requires an understanding of unbounded
$\WOT$-convergence. For example,
suppose $\I$ belongs to $\Id_{l}(\gL_{n})$ with $\mu\,(\I)=\H_{n}$. It is not
even known if one must have $\I=\gL_{n}$. For, one would like to say that the
identity $I$
belongs to  $\I$, but all that can be said is $\xi_{1}=\lim_{\alpha}J_{\alpha}
\xi_{1}$
for some $J_{\alpha}$ in $\I$.
 For bounded nets, $\WOT$-convergence amounts to strong convergence on the
vector $\xi_{1}$ \cite{topalg}.
However, this is not true for unbounded nets. Indeed, as an example consider
the sequence $J_{m}$ of operators in $\gL_{n}$ given by
\[ J_{m} = \sum_{k=0}^{m}\frac{1}{k+1}L_{1}^{m}. \]
It is clear that
\[ \lim_{m\rightarrow\infty}J_{m}\xi_{1} = \sum_{k\geq 0}\frac{1}{k+1}\xi_{
z_{1}^{k}}, \]
but the latter vector does not represent the Fourier coefficients of any
operator in $\gL_{n}$.
\end{rem}

Nonetheless, it has been shown that the maps $\mu$ and $i$ define a bijective
correspondence between the cyclic subspaces of $\Lat (\gL_{n})$ determined by
cyclic flips on the one hand, and the image under $i$ in ${\rm Id}_{l}(\gL_{n})$ of these subspaces on the other. Concerning the lattice properties of these
maps, it is not hard to show that $\mu$ sends closed spans to $\WOT$-closed
sums and $i$ sends intersections to intersections. However, the behaviour of
$\mu$ on intersections and $i$ on sums again comes back to requiring an
understanding of unbounded $\WOT$-convergence.

\end {document}